\documentclass[11pt]{article}

\usepackage{amsmath}
\usepackage{amssymb}
\usepackage{eucal}

\begin{document}

\baselineskip 16pt

\title{ On $\sigma$-arithmetic  graphs of finite groups
 }

\author{ \\
{ Alexander  N. Skiba}\\
{\small Department of Mathematics and Technologies of Programming,  Francisk Skorina Gomel State University,}\\
{\small Gomel 246019, Belarus}\\
{\small E-mail: alexander.skiba49@gmail.com}}

\date{}
\maketitle

\begin{abstract} Let   $G$ be a finite group and $\sigma$ a partition of the set of all
 primes  
$\Bbb{P}$, that is,  $\sigma =\{\sigma_{i} \mid
 i\in I \}$, where   $\Bbb{P}=\bigcup_{i\in I} \sigma_{i}$
 and $\sigma_{i}\cap \sigma_{j}= \emptyset  $ for
 all $i\ne j$.  If $n$ is an integer, we write $\sigma (n)=\{\sigma_{i} \mid \sigma_{i}\cap \pi (n)\ne \emptyset  \}$  and
   $\sigma (G)=\sigma (|G|)$.

  We say that a 
chief factor $H/K$ is: \emph{$\pi$-central} in $G$
 ($\pi \subseteq \Bbb{P}$) if $H/K$ and $G/C_{G}(H/K)$ are 
$\pi$-groups; $\sigma$-central in $G$ if 
 $H/K$ is  $\sigma _{i}$-central in $G$ for some $i$; the symbol $F_{\{\sigma_{i} \}}(G)$
 denotes the largest normal $\{\sigma_{i} \}$-nilpotent subgroup of $G$.

The group $G$ is said to be: 
$\sigma$-soluble  if every 
chief factor $H/K$ of $G$ is a $\sigma _{i}$-group for some $i=i(H/K)$; 
 \emph{$\sigma$-nilpotent} if every chief factor of 
$G$ is $\sigma$-central in $G$; 
     $\{\sigma _{i}\}$-nilpotent  if  every  chief factor $H/K$ of $G$ with
 $\sigma (H/K)\cap \sigma _{i}\ne \emptyset$
 is $\sigma$-central in $G$.

 The symbol 
${\mathfrak{N}}_{\sigma }$ denotes the class of all $\sigma$-nilpotent 
groups; an
 \emph{${\mathfrak{N}}_{\sigma }$-critical group } is a non-$\sigma$-nilpotent group
 whose all proper subgroups are $\sigma$-nilpotent.  
  
We call
 a graph $\Gamma$ with the set of all vertices $V(\Gamma)=\sigma (G)$ ($G\ne 1$) a
 \emph{$\sigma$-arithmetic  graph} of $G$, and we associate with $G\ne 1$ the
 following three directed $\sigma$-arithmetic graphs:  (1) the  \emph{$\sigma$-Hawkes graph} $\Gamma _{H\sigma }(G)$ of $G$ 
is a $\sigma$-arithmetic  graph of $G$ in which $(\sigma _{i}, \sigma 
_{j})\in E(\Gamma _{H\sigma }(G))$ if $\sigma _{j}\in \sigma 
(G/F_{\{\sigma _{i}\}}(G))$; (2) the  \emph{$\sigma$-Hall graph} 
 $\Gamma _{{\sigma }Hal}(G)$ of $G$ in which  
$(\sigma _{i}, \sigma 
_{j})\in E(\Gamma _{{\sigma }Hal}(G))$ if  for some Hall 
$\sigma _{i}$-subgroup $H$ of $G$ we have  $\sigma _{j}\in \sigma 
(N_{G}(H)/HC_{G}(H))$; 
(3)  the  \emph{$\sigma$-Vasil'ev-Murashko graph}   $\Gamma _{{\mathfrak{N}_{\sigma }}}(G)$
 of $G$ in which  
$(\sigma _{i}, \sigma 
_{j})\in E(\Gamma _{{\mathfrak{N}_{\sigma }}}(G))$ if 
for some  ${\mathfrak{N}_{\sigma }}$-critical subgroup $H$ 
of $G$ we have $\sigma _{i} \in \sigma (H)$ and  
$\sigma _{j}\in \sigma (H/F_{\{\sigma _{i}\}}(H))$.

In this paper, we study the structure of
 $G$ depending on the properties of these three  graphs of $G$. 
 In particular, we prove the following result.

 {\bf Theorem 1.4.}  {\sl  Let $G\ne 1$. Then the following statements are 
equivalent:}

(1) {\sl $G$ has a normal series $1= G_{1} <  G_{2}  < \cdots   <  G_{t-1} < 
G_{t}=G$ in which for every  $i=1, \ldots , t$, $G_{i}/G_{i-1}$ is a $\sigma _{i_{j}}$-group 
and $G_{i-1}$ and $G/G_{i}$ are $\sigma _{i_{j}}'$-groups for some
 $\sigma _{i_{j}}\in \sigma (G)$;}

(2) {\sl $\Gamma _{H\sigma }(G)$ has no circuits;}

(3) {\sl $G$ is $\sigma$-soluble and  $\Gamma _{{\mathfrak{N}}_{\sigma }}(G)$ 
 has no circuits. }

\end{abstract}

\let\thefootnoteorig\thefootnote
\renewcommand{\thefootnote}{\empty}


\let\thefootnote\thefootnoteorig

\section{Introduction}

Throughout this paper, all groups are finite and $G$ always denotes
a finite group. Moreover,    $\mathbb{P}$ is the set of all  primes,
 $\pi \subseteq \mathbb{P} $ and $\pi'=\mathbb{P}\setminus \pi$. If 
 $n$ is an integer, then the symbol $\pi (n)$ denotes
 the set of all primes dividing $n$; as usual, $\pi (G)=\pi (|G|)$, the set of all
  primes dividing the order of the group $G$.   We say that a 
chief factor $H/K$ is \emph{$\pi$-central} (in $G$) if $H/K$ and $G/C_{G}(H/K)$ are 
$\pi$-groups.

We write $V(\Gamma)$ to denote the set of all vertices of the graph $\Gamma$, 
$E(\Gamma)$ is the set of all edges of   $\Gamma$. If $\Gamma _{1} $ is a graph such that 
$V(\Gamma _{1})\subseteq  V(\Gamma)$ and  $E(\Gamma)\subseteq E(\Gamma 
_{1})$, then $\Gamma _{1}$ is a \emph{subgraph} of  $\Gamma$ and 
in this case we write $\Gamma _{1}\subseteq  \Gamma$.  
If $\Gamma _{1}$ and $\Gamma _{2}$ are graphs, then  their \emph{union}
$\Gamma =\Gamma _{1} \bigcup \Gamma _{2}$ is a  graph such 
that $V(\Gamma)=V(\Gamma _{1}) \bigcup V(\Gamma _{2})$ and 
$E(\Gamma)=E(\Gamma _{1}) \bigcup E(\Gamma _{2})$.
 
A graph $\Gamma$ with $V(\Gamma)=\pi (G)$ is  called an \emph{arithmetic graph of $G$}
\cite{mur-vas}. 
From the properties of arithmetic graphs of the group $G$, you can extract 
meaningful information about the structure of $G$ (see Introductions in \cite{mur-vas, mur}).

A classic example of an arithmetic graph of  $G$ is the Gr\"{u}enberg-Kegel 
graph $\Gamma _{p}(G)$ in which $(q, r)\in E(\Gamma _{p}(G))$ if $G$ has
 an element of order $qr$.
Graphs of this type proved useful in analyzing 
 of many questions of 
the group theory (see, in particular,  \cite{will, kon, maz}). The    
\emph{Hawkes graph} $\Gamma _{H}(G)$ of $G$ is another example 
arithmetic graph of $G$ in which $(q, r)\in E(\Gamma _{H}(G))$ if 
$q\in \pi (G/F_{p}(G))$, where $F_{p}(G)=O_{p', p}(G)$ is the largest normal $p$-nilpotent subgroup of
 $G$. Hawkes proved \cite{Haw}  that $G$ is a Sylow tower group if and 
only if the graph $\Gamma _{H}(G)$ has no circuits. Further applications 
of such a graph were found in the recent papers \cite{mur-vas, mur}. Note 
that  the   Hawkes graph of $G$  is  directed.  Two another interesting 
 directed arithmetic  graphs of a group
 were used in the papers  \cite{Kas, al, mur-vas, mur}. First recall that 
the 
\emph{Sylow graph} $\Gamma_{S} (G)$ of $G$ \cite{Kas} is an arithmetic graph of $G$ in 
which $(p, q)\in E(\Gamma_{S}(G)) $ if $q\in \pi (N_{G}(P)/PC_{H}(P))$ for some Sylow $p$-subgroup 
$P$ of $G$. Kazarin and others proved in \cite{Kas} that $\Gamma_{S}(G)$ is 
connected and its diameter is at most  5 whenever $G$ is an almost simple group. 
Recall also that the \emph{$\mathfrak{N}$-critical } or the 
\emph{Vasil'ev-Muraschko graph } $\Gamma_{\mathfrak{N}c}(G)$ of $G$  \cite{mur-vas}  is
 an arithmetic graph of $G$ in 
which $(p, q)\in E(\Gamma_{\mathfrak{N}c}(G))$ if for some $p$-closed  Schmidt subgroup $A$ 
of $G$  we have $\pi (A)=\{p, q\}$. Vasil'ev and Murashko proved in \cite{mur-vas} that if 
$\pi _{1}, \ldots , \pi _{n}$ are sets of vertices of connected components 
in 
$\Gamma_{\mathfrak{N}c}(G)$, then $G=O_{\pi _{1}}(G) \times \cdots \times O_{\pi _{n}}(G)$.

Now let $\sigma$  any partition of  
$\Bbb{P}$, that is,  $\sigma =\{\sigma_{i} \mid
 i\in I \}$, where   $\Bbb{P}=\bigcup_{i\in I} \sigma_{i}$
 and $\sigma_{i}\cap \sigma_{j}= \emptyset  $ for
 all $i\ne j$.  
When studying the $\sigma$-properties of a group $G$, that are, properties 
of $G$ which are defined by the choice of $\sigma$,  it is convenient to 
use the $\sigma$-analogues of arithmetic graphs of $G$. 

By analogy with the notations $\pi (n)$  and $\pi (G)$,
 we write $\sigma (n)=\{\sigma_{i} \mid \sigma_{i}\cap \pi (n)\ne \emptyset  \}$  and
   $\sigma (G)=\sigma (|G|)$ \cite{alg12}. We call any graph $\Gamma$ with $V(\Gamma)=\sigma (G)$ ($G\ne 1$) a
 \emph{$\sigma$-arithmetic  graph of $G$}.

Before continuing, we recall some basic concepts and notations of  the theory of 
$\sigma$-properties of  groups.

The group $G$ is said to be \cite{1}:  
\emph{$\sigma$-soluble}  if every 
chief factor $H/K$ of $G$ is a $\sigma _{i}$-group for some $i=i(H/K)$; 
 \emph{$\sigma$-nilpotent} if every chief factor $H/K$ of 
$G$ is $\sigma$-central in $G$, that is, 
 $H/K$ is  $\sigma _{i}$-central in $G$ for some $i=i(H/K)$; 
     \emph{$\{\sigma _{i}\}$-nilpotent}  if  every  chief factor $H/K$ of $G$ with
 $\sigma (H/K)\cap \sigma _{i}\ne \emptyset$  is $\sigma$-central in $G$ ; the symbol $F_{\{\sigma_{i} \}}(G)$
 denotes the largest normal $\{\sigma_{i} \}$-nilpotent subgroup of $G$.

 The symbol 
${\mathfrak{N}}_{\sigma }$ denotes the class of all $\sigma$-nilpotent 
groups; an
 \emph{${\mathfrak{N}}_{\sigma }$-critical group } is a non-$\sigma$-nilpotent group
 whose all proper subgroups are $\sigma$-nilpotent. A \emph{Schmidt group
} is a non-nilpotent group  whose all proper subgroups are 
nilpotent.

Now we associate with $G\ne 1$ three directed $\sigma$-arithmetic graphs.

{\bf Definition 1.1.} (1) The  \emph{$\sigma$-Hawkes graph} $\Gamma _{H\sigma }(G)$ of $G$ 
is a $\sigma$-arithmetic  graph of $G$ in which $(\sigma _{i}, \sigma 
_{j})\in E(\Gamma _{H\sigma }(G))$ if $\sigma _{j}\in \sigma (G/F_{\{\sigma _{i}\}}(G))$.

(2) The  \emph{$\sigma$-Hall graph}  $\Gamma _{{\sigma }Hal}(G)$ of $G$ in which  
$(\sigma _{i}, \sigma 
_{j})\in E(\Gamma _{{\sigma }Hal}(G))$ if  for some Hall 
$\sigma _{i}$-subgroup $H$ of $G$ we have  $\sigma _{j}\in \sigma 
(N_{G}(H)/HC_{G}(H))$. 
                           
(3)  The  \emph{$\sigma$-Vasil'ev-Murashko graph}   $\Gamma _{{\mathfrak{N}_{\sigma }}}(G)$
 of $G$ in which  
$(\sigma _{i}, \sigma 
_{j})\in E(\Gamma _{{\mathfrak{N}_{\sigma }}}(G))$ if 
for some  ${\mathfrak{N}_{\sigma }}$-critical subgroup $H$ 
of $G$ we have $\sigma _{i} \in \sigma (H)$ and  
$\sigma _{j}\in \sigma (H/F_{\{\sigma _{i}\}}(H))$.

Note that from the properties of Schmidt groups \cite[III, Satz 5.2]{hupp} it follows   in the case when
 $\sigma =\sigma ^{1}=\{\{2\}, \{3\}, \ldots \}$
 (we use here the notations in \cite{alg12}) we have $\Gamma _{H\sigma }(G)=\Gamma _{H}(G)$, 
$\Gamma _{{\sigma }Hal}(G)=\Gamma _{S}(G)$ and
  $\Gamma _{{\mathfrak{N}_{\sigma }}}(G)=\Gamma _{\mathfrak{N}c}(G)$.

{\bf Proposition 1.2.} {\sl  Let $G \ne 1$. Then }

(1) {\sl $$\Gamma _{{\sigma }Hal}(G)  \subseteq
 \Gamma _{{\mathfrak{N}}_{\sigma }}(G)\subseteq   \Gamma _{H\sigma 
}(G)\eqno(*).$$ Moreover, 
in general, both inclusions   in (*) may be strict, and }

(2) {\sl the equalities $$\Gamma _{{\sigma 
}\text{Hal}}(G) =
 \Gamma _{{\mathfrak{N}}_{\sigma }}(G)=  \Gamma _{H\sigma 
}(G)$$ hold if and only if the graph $\Gamma _{H\sigma 
}(G)$ has no loops, that is, $G$ is $\sigma$-soluble  with  $l_{\sigma 
_{i}}(G)\leq 1$ } (see \cite[p. 249]{rob}) {\sl for all  $\sigma _{i}\in \sigma (G)$}.

Proposition 1.2 is a  motivation for the following 

{\bf Question 1.3.} What is the structure of the group $G\ne 1$ in which the 
equality $\Gamma _{{\sigma }Hal}(G)=\Gamma _{{\mathfrak{N}}_{\sigma }}(G)$ holds?

Our first goal is to prove the following analogous of the above-mentioned 
Hawkes result.
   
 {\bf Theorem 1.4.}  {\sl  Let $G\ne 1$. Then the following statements are 
equivalent:}

(1) {\sl $G$ has a normal series $1= G_{1} <  G_{2}  < \cdots   <  G_{t-1} < 
G_{t}=G$ in which for every  $i=1, \ldots , t$, $G_{i}/G_{i-1}$ is a $\sigma _{i_{j}}$-group 
and $G_{i-1}$ and $G/G_{i}$ are $\sigma _{i_{j}}'$-groups for some
 $\sigma _{i_{j}}\in \sigma (G)$;}

(2) {\sl $\Gamma _{H\sigma }(G)$ has no circuits;}

(3) {\sl $G$ is $\sigma$-soluble and  $\Gamma _{{\mathfrak{N}}_{\sigma }}(G)$ 
 has no circuits. }

 From Proposition 1.2(2) and Theorem 1.4 we get

{\bf Corollary 1.5} (Hawkes \cite{Haw}). {\sl  Let $G\ne 1$. Then 
 $G$ is a Sylow tower group if and 
only if the graph $\Gamma _{H}(G)$ has no circuits. }

{\bf Corollary 1.6} (See Theorem 7.1 in \cite{mur-vas}). {\sl  Let $G\ne 1$. 
If $\Gamma _{\mathfrak{N}c}(G)$ has no circuits, then $G$ is a Sylow tower group}.

The integers $n$ and $m$ are called \emph{$\sigma$-coprime} if $\sigma (n) \cap 
\sigma (m)=\emptyset. $

 {\bf Theorem 1.7.}  {\sl  Let $1 \ne G=AB=BC=AC$, where $A$, $B$ and $C$ are $\sigma$-soluble 
subgroups of $G$.  Then the following statements hold.}

(1) {\sl If $G$ is $\sigma$-soluble, then
 $\Gamma _{{\mathfrak{N}_{\sigma }}}(A)\bigcup\Gamma _{{\mathfrak{N}_{\sigma }}}(B) 
 \bigcup \Gamma _{{\mathfrak{N}_{\sigma }}}(C)=  \Gamma _{{\mathfrak{N}_{\sigma }}}(G).$ }

(2) {\sl If the  indices $|G:A|, |G:B|,
 |G:C|$ are pair   $\sigma$-coprime, then $\Gamma _{H\sigma }(A) \bigcup\Gamma _{H\sigma }(B)
  \bigcup \Gamma _{H\sigma }(C)= \Gamma _{H\sigma }(G).$}

Since every group $G=AB=BC=AC$, where $A$, $B$ and $C$ are soluble 
subgroups of $G$, is soluble by \cite{kazarin}, we get from Theorem 1.7 the following

 {\bf Corollary 1.8} (See Theorem 6.2 in \cite{mur-vas})  {\sl
  Let $1 \ne G=AB=BC=AC$, where $A$, $B$ and $C$ are soluble 
subgroups of $G$.  Then:}

(1) {\sl 
 $\Gamma _{{\mathfrak{N}c}}(A)\bigcup\Gamma _{{\mathfrak{N}c}}(B) 
 \bigcup \Gamma _{{\mathfrak{N}c}}(C)=  \Gamma _{{\mathfrak{N}c}}(G),$  and }

(2) {\sl if the  indices $|G:A_{1}|, |G:A_{2}|,
 |G:A_{3}|$ are pair   coprime, then $\Gamma _{H}(A) \bigcup\Gamma _{H}(B)
  \bigcup \Gamma _{H}(C)= \Gamma _{H}(G).$}

In what follows, $\Pi \subseteq \sigma$ and $\Pi'= \sigma \setminus \Pi$. 
An integer $n$ is called a \emph{$\Pi$-number} if $\sigma (n)\subseteq \Pi$. A 
subgroup $H$ of $G$  is said to be  a \emph{Hall $\Pi$-subgroup} of $G$ if
 $|G:H|$ is a $\Pi'$-number and 
 $H$ is a \emph{$\Pi$-group},
 that is, $|H|$ 
is a $\Pi$-number; $G$ is called \emph{$\Pi$-closed} if $G$ has a 
normal Hall $\Pi$-subgroup.  

 In proving of Theorems 1.4 and 1.7, the next two propositions turned out to be useful.

{\bf Proposition 1.9.}  {\sl  Let $G\ne 1$ and $\Pi_{1}=\Pi\cap \sigma 
(G)$ and  $\Pi_{2}=\Pi' \cap \sigma (G)$.  Let $\Gamma$ be a  
$\sigma$-arithmetic   
 graph of $G$. Suppose that for every $\sigma _{i}\in \Pi_{2}$ 
there is no $\sigma _{j}\in \Pi_{1}$ such that $(\sigma _{i}, \sigma _{j})\in E(\Gamma)$.  
 If either  $\Gamma=\Gamma _{H\sigma }(G)$
 or $G$ is $\sigma$-soluble and   $\Gamma=\Gamma _{{\mathfrak{N}}_{\sigma 
}}(G)$, then $G$ is $\Pi _{1}$-closed. }

{\bf Corollary 1.10} (See Theorem 5.1 in \cite{mur-vas}).  {\sl
  Let $G\ne 1$ and $\pi_{1}=\pi\cap \pi 
(G)$ and  $\pi_{2}=\pi' \cap \pi (G)$.   Suppose that for every $p\in \pi_{2}$ 
there is no $p\in \pi_{1}$ such that $(p, q)\in E(\Gamma _{H}(G))$.  
 Then $G$ is $\pi _{1}$-closed. }

{\bf Proposition 1.11.}  {\sl Suppose that $G$ is a $\sigma$-soluble 
minimal non-$\Pi$-closed 
group, that is, $G$  is not $\Pi$-closed but every proper subgroup of $G$ is 
$\Pi$-closed. Then $G$ is a $\Pi'$-closed  Schmidt group. }

Finally, we prove the following

{\bf Theorem 1.12.}  {\sl  Let $G\ne 1$.  Then the following statements are equivalent:   }

(1) {\sl $G$ is $\sigma$-nilpotent;}

(2) {\sl  Each vertex of the graph $\Gamma _{H\sigma }(G)$ is isolated.    
 }

(3) {\sl  Each vertex of the graph $\Gamma _{{\mathfrak{N}}_{\sigma }}(G)$ 
is isolated. }

 (4) {\sl  $G$ is $\sigma$-soluble and 
 each vertex of the graph $\Gamma _{{\sigma }Hall}(G)$ is 
isolated.}

{\bf Corollary  1.13}  (See Theorem 5.4 in \cite{mur-vas}).  {\sl  Let $G\ne 1$. 
 If 
$\pi _{1}, \ldots , \pi _{n}$ are sets of vertices of connected components 
in 
$\Gamma_{\mathfrak{N}c}(G)$, then $G=O_{\pi _{1}}(G) \times \cdots \times 
O_{\pi _{n}}(G)$. }

{\bf Proof.}  Let
 $\sigma =\{\pi _{1}, \ldots , \pi _{n}, \pi _{0}'\}$,
 where $\pi_{0}=\pi _{1}\cup \cdots \cup \pi 
_{n}$. Then from the hypothesis and Theorem  1.12 it follows that $G$ is $\sigma$-nilpotent, 
that is, $G=O_{\pi _{1}}(G) \times \cdots \times 
O_{\pi _{n}}(G)$.  The corollary is proved.

 \section{$\{\pi\}$-nilpotent and $\sigma$-nilpotent groups}

We say that $G$ is \emph{$\{\pi\}$-nilpotent} if every chief factor $H/K$ 
of $G$ such that $\pi \cap \pi(H/K)\ne \emptyset$ is $\pi$-central in $G$.

In what follows,  $\mathfrak{F}$ is a  class of groups containing all identity groups;  
  $G^{\frak {F}}$  denotes the intersection of all normal
subgroups $N$ of $G$ with  $G/N\in {\frak {F}}$; $G_{\frak {F}}$  is the product
  of all normal
subgroups $N$ of $G$ with  $N\in {\frak {F}}$.  The class $\mathfrak{F}$ is  said to be:
 a \emph{formation}  
if for every group $G$ every homomorphic image of $G/G^{\frak {F}}$
 belongs to $\frak {F}$;  a \emph{Fitting class}  if  
  for every group $G$   every normal subgroup of $G_{\frak {F}}$
 belongs to $\frak {F}$.  
 The class   $\mathfrak{F}$
 is called: \emph{saturated}  if 
 $G\in \mathfrak{F}$ whenever $G/\Phi (G)\in \mathfrak{F}$; \emph{hereditary} (Mal'cev \cite{mal}) if
 $H\in \frak {F}$ whenever $H\leq G \in
\frak {F}$.  
  
{\bf Proposition 2.1.} {\sl The class of all   $\{\pi\}$-nilpotent groups
 ${\mathfrak{G}}_{\{\pi\}}$ is a hereditary  Fitting formation.}

{\bf Proof.}  First note 
that the class   ${\mathfrak{G}}_{\{\pi\}}$ is obviously closed under taking
homomorphic images.

Next we show that 
${\mathfrak{G}}_{\{ \pi\}}$  is a hereditary  Fitting class.
 Let $E\leq G$ and  $R$ and $N$  normal subgroups  of $G$. 
First assume that $G\in  {\mathfrak{G}}_{\{\pi\}}$ and let $1=G_{0} < G_{1} < \cdots < 
G_{t-1} < G_{t}=G$ be a chief series of $G$. Consider the series  
$$1=G_{0}\cap E < G_{1}\cap E < \cdots < 
G_{t-1}\cap E < G_{t}\cap E=E$$ in $E$ and let $H/K$ be a chief factor of $E$ such 
that  $\pi\cap \pi (H/K)\ne \emptyset$ and   $ G_{k-1}\cap E \leq  
 K  <   H  \leq G_{k}\cap E $ for some $k$. Then $\pi\cap \pi ((G_{k}\cap E)/(G_{k-1}\cap E))\ne 
\emptyset$ and from 
the isomorphism $$(G_{k}\cap E)/(G_{k-1}\cap E) \simeq  G_{k-1}(G_{k}\cap E)/G_{k-1}=
    G_{k}/G_{k-1}$$ we  get that  $\pi\cap \pi (G_{k}/G_{k-1})\ne \emptyset$,  
     so $G_{k}/G_{k-1}$ and  $G/C_{G}(G_{k}/G_{k-1})$ are $\pi$-groups by 
hypothesis. Then $H/K$
 and $$E/(E\cap C_{G}(G_{k}/G_{k-1}))\simeq EC_{G}(G_{k}/G_{k-1})/C_{G}(G_{k}/G_{k-1})$$ are 
 $\pi$-groups.  
 But $$E\cap C_{G}(G_{k}/G_{k-1
})=E\cap C_{G}((G_{k}\cap E)/(G_{k-1
}\cap E))   \leq 
C_{E}(H/K)$$ and so $H/K$ is $\pi$-central in $E$. Hence every  chief 
factor $T/L$ of $E$ such that $\sigma _{i}\cap \pi (T/L)\ne \emptyset$ is $\pi$-central in $E$
 by the Jordan-H\"{o}lder theorem for 
the chief series, that is, $E$ is  $\{\pi\}$-nilpotent. 
 Therefore the class 
${\mathfrak{G}}_{\{\pi\}}$  is  hereditary.

Now we show  that  if $R, N \in {\mathfrak{G}}_{\{\pi\}}$, then $R N
 \in {\mathfrak{G}}_{\{ \pi\}}$. We can assume without loss of 
generality that $RN=G$. Let $H/K$ be any chief factor of 
$G$  such 
that  $\pi\cap \pi (H/K)\ne \emptyset$. In view of the Jordan-H\"{o}lder theorem for 
the chief series, $G$ has a chief factor $T/L$ such that $H/K$ and $T/L$ are  
$G$-isomorphic and either $T\leq N$ or $N\leq L$.

 First we show  that if 
 $H\leq R$, then $H/K$ and $ RC_{G}(H/K)/C_{G}(H/K)$ are $\pi$-groups. 
   Indeed,      $H/K=(H_{1}/K)\times \cdots \times 
(H_{t}/K)$, where $H_{  {i}}/K$ is a minimal normal subgroup of $R/K$ for all 
$i=1, \ldots, t$, by \cite[Chapter A, Proposition 4.13(c)]{DH} and  so   
$$C_{R} (H_{1}/K)\cap  \cdots \cap C_{R}(H_{t}/K) =   C_{R}(H/K)=C_{G}(H/K)\cap R.$$  
 But $R/K$ is $\{\pi\}$-nilpotent  since  $R \in {\mathfrak{G}}_{\{ \pi\}}$ and 
  the class   ${\mathfrak{G}}_{\{\pi\}}$ is closed under taking
homomorphic images. 
 Hence $$R/(C_{R} (H_{1}/K)\cap  \cdots \cap C_{R}(H_{t}/K))=R/(R\cap C_{G}(H/K))\simeq 
RC_{G}(H/K)/C_{G}(H/K)$$ is a $\pi$-group.  Since $R$ is  $\{\pi\}$-nilpotent and
 $\pi_{i}\cap \pi (H/K)\ne \emptyset$, it follows that $H/K$
  is a $\pi$-group.
 Similarly, if 
 $ T\leq N$, then $H/K$ and $$NC_{G}(H/K)/C_{G}(H/K)= NC_{G}(T/L)/C_{G}(T/L)$$ are
 $\pi$-groups and
 so   
 $$G/C_{G}(H/K)=(R/C_{G}(H/K)/C_{G}(H/K))(N/C_{G}(H/K)/C_{G}(H/K))$$  is a 
$\pi$-group.  Hence in this case $H/K$ is $\pi$-central in $G$.

Now  assume that $R\leq K$ or $N\leq L$, $R\leq K$ say.
 Then $H=R(H\cap N)=K(H\cap N)$ and so from the 
$G$-isomorphism   $H/K\simeq (H\cap N)/(K\cap N)$  we get that  $H/K$ is a 
$\pi$-group and that 
$C_{G}(H/K)=C_{G}((H\cap N)/(K\cap N))$, where 
$$NC_{G}((H\cap N)/(K\cap N))/C_{G}((H\cap N)/(K\cap N))$$ is  a
 $\pi$-group by the previous paragraph.  On the other hand, $R\leq C_{G}(H/K)$. Hence 
$$RNC_{G}((H\cap N)/(K\cap N))/C_{G}((H\cap N)/(K\cap N))   =G/C_{G}((H\cap N)/(K\cap N))=  
  G/C_{G}(H/K)$$   is  a $\pi$-group. Therefore $H/K$ is $\pi$-central in $G$.  
   Hence  $RN\in  {\mathfrak{G}}_{\{\pi\}}$, so 
  ${\mathfrak{G}}_{\{\pi\}}$  is a  Fitting  class.

Now we have only to show that if
 $G/R, G/N \in {\mathfrak{G}}_{\{\pi\}}$ and  $H/K$
 is  a chief factor   of $G$ such 
that  $\pi\cap \pi (H/K)\ne \emptyset$ and 
 $R\cap N\leq K$, then $H/K$ is $\pi$-central in $G$.  
If $RK\ne RH$, then from the $G$-isomorphisms $$HR/KR\simeq H/(H\cap KR)=H/K(H\cap R)=H/K$$  and 
$(HR/R)/(KR/R)\simeq HR/KR$  we get that $H/K$ is $\pi$-central in $G$.   Finally, 
assume that  $RK= RH$. Then from the  $G$-isomorphisms $H/K\simeq (H\cap R)/(K\cap R)$ and
  $ (H\cap R)N/(K\cap R)N\simeq (H\cap R)/(K\cap R)$ we again get that  
$H/K$ is $\pi$-central in $G$. 
The proposition is proved.

{\bf Corollary 2.2.} {\sl The class of all   $\sigma$-nilpotent groups
 ${\mathfrak{N}}_{\sigma}$ is a hereditary  Fitting formation.}

{\bf Proof. } This follows from Proposition 2.1 and the fact that
 ${\mathfrak{N}}_{\sigma}=\bigcap
 _{i \in I}{\mathfrak{G}}_{\{\sigma _{i}\}}$.

We say that a maximal subgroup $M$ of $G$ is \emph{$\pi$-normal in $G$} if either 
$|G:M|$ is a $\pi'$-number or some  chief factor $V/M_{G}$ of $G$ is $\pi$-central in $G$. 

{\bf Proposition 2.3.} {\sl The following conditions are equivalent:}

(1) {\sl $G$ is $\{\pi\}$-nilpotent}

(2) {\sl  $G$ has a normal $\pi$-complement.}

(3) {\sl  $G$ is $\pi$-separable and every maximal subgroup of $G$ is $\pi$-central in $G$. }

{\bf Proof. }  (1) $\Rightarrow$  (2)  Assume that this implication is false and let 
$G$ be a counterexample of minimal order. Then $\pi \cap   \pi (G)   \ne \emptyset \ne 
\pi' \cap   \pi (G)$. 

Let $R$ be a minimal normal 
subgroup of $G$ and $C=C_{G}(R)$. Then $G/R$ is  $\{\pi\}$-nilpotent by Proposition 2.1,
 so $G/R$ has a normal 
$\pi$-complement $V/R$, that is, a normal  Hall $\pi'$-subgroup of $G/R$. If $R$ is
 a $\pi'$-group, then $V$ is a  Hall $\pi'$-subgroup of $G$. Hence we may assume that
 every minimal normal 
subgroup of $G$ is a $\pi$-group  since $G$ is $\pi$-separable. 
Therefore $G/C$ is also a $\pi$-group. Assume that $C\ne G$. In view of 
Proposition 2.1, $C$ is $\{\pi\}$-nilpotent and so $C$  has a normal 
$\sigma _{i}$-complement $U$ by the choice of $G$. But $U$ is characteristic in $C$ and  
$U$ is a Hall $\pi'$-subgroup of $G$. Therefore $U$ is a normal 
$\pi$-complement of $G$, contrary to our assumption on $G$. 
Therefore $C=G$, so $R\leq Z(G)\cap V\leq Z(V)$. Let $E$ be a Hall $\pi'$-subgroup of $V$. Then $V=R\times E$, so $E$ is a a normal 
$\pi$-complement of $G$. This contradiction completes the proof of 
the implication (1) $\Rightarrow$  (2)
 
 (2) $\Rightarrow$  (1), (3)  Assume that  $G$ has a normal 
$\pi$-complement $V$. Then for every chief factor $H/K$ of $G$ 
such that $V\leq K$ we have $V\leq C_{G}(H/K)$, so $H/K$ and  
$G/C_{G}(H/K)$ are $\pi$-group. On the other hand, for every chief factor $H/K$ of $G$ 
such that $H\leq V$ we have $\pi (H/K)\cap \pi=\emptyset$, so  
every chief factor $H/K$ of $G$ such that $\pi (H/K)\cap \pi\ne \emptyset$ i
s $\pi$-central in $G$ by the Jordan-H\"{o}lder theorem for the chief series.
 Hence $G$ is $\{\pi\}$-nilpotent, so (2) $\Rightarrow$  (1). Moreover, 
if $M$ is a maximal subgroup of $G$ such that $|G:M|$ is not 
$\pi'$-number, then for every chief factor $V/M_{G}$ we have $G=MV$ and so 
from $|G:M|=|MV|:|M|=|V:V\cap M|$  we get that $\pi(V/M_{G})\cap \pi\ne \emptyset$. Hence  
$V/M_{G}$ of $G$ is $\pi$-central in $G$. It is  also clea that $G$ is 
$\pi$-separable.  Hence  (2) $\Rightarrow$  (3).

 (3) $\Rightarrow$  (2)  Assume that this implication is false and let 
$G$ be a counterexample of minimal order.
Let $R$ be a minimal normal  subgroup of $G$ and $C=C_{G}(R)$. Then $R$ is either a 
$\pi$-group or a $\pi'$-group. Moreover, the hypothesis holds for 
$G/R$, so  $G/R$ has a normal 
$\pi$-complement $V/R$ by the choice of $G$. If $R$ is a $\pi'$-group, then $V$ is a 
normal $\pi$-complement in $G$. Therefore   $R$ is a $\pi$-group. Then $R$ has a complement 
$E$ in $V$ and every two such complements are conjugate in $V$
 by the Schur-Zassenhaus theorem, so $G=VN_{G}(E) = REN_{G}(E) =RN_{G}(E)$ by the Frattini Argumant. 
 Hence  $R\nleq  \Phi (G)$. If $R/1 $ is $\pi$-central in $G$, 
then  $E\leq C_{G}(R)$ and hence $V=R\times E$. But then $E$ is a normal 
$\pi$-complement in $G$. Therefore   $R/1 $ is not $\pi$-central in $G$.

Let $M$ be a maximal subgroup of $G$ such that $G=RM$. 
Then $|G:M|$ divides $|R|$, so some 
chief factor $V/M_{G}$ of $G$ is $\pi$-central in $G$, that is, $V/M_{G}$ and
 $G/C_{G}(V/M_{G})$ are $\pi$-groups. Then from the $G$-isomorphism 
$R\simeq RM_{G/}M_{G}$ we get that $RM_{G}/M_{G}\ne V/M_{G}$. Hence 
$C_{G/M_{G}}(V/M_{G})=RM_{G}/M_{G}$ by \cite[Chapter A, Theorem 15.2(3)]{DH}, so 
$(G/M_{G})/(RM_{G}/M_{G})$ and hence $G/M_{G}$ are $\pi$-groups. But then 
$RM_{G/}M_{G}$ and so $R/1$ are $\pi$-central in $G$. This contradiction 
completes the proof of the implication.

The proposition is proved.

{\bf Corollary 2.4.} {\sl The classes ${\mathfrak{G}}_{\{\pi\}}$ and
 ${\mathfrak{N}}_{\sigma}$  are  saturated.}

{\bf Proof. }  First assume that $G/\Phi (G)\in {\mathfrak{G}}_{\{\pi\}}$. 
Then $G/\Phi (G)$ is  $\pi$-separable and every maximal subgroup of $G/\Phi (G)$ is 
$\pi$-central in $G/\Phi (G)$ by Proposition 2.3. It follows that  $G$ is $\pi$-separable and every maximal 
subgroup of $G$ is $\pi$-central in $G$. Therefore  $G\in 
{\mathfrak{G}}_{\{\pi\}}$ by Proposition 2.3. Hence ${\mathfrak{G}}_{\{\pi\}}$ is saturated.            But then
 ${\mathfrak{N}}_{\sigma}=\bigcap
 _{i \in I}{\mathfrak{G}}_{\{\sigma _{i}\}}$ is saturated.  The lemma is 
proved.

\section{Proofs of Propositions 1.2, 1.9 and 1.11}

{\bf Lemma 3.1} (See  \cite[III, 5.2]{hupp}).
  {\sl If  
$G$ is a   
 Schmidt group, then $G=P\rtimes Q$, where    $P=G^{\frak{N}}$ 
 is a Sylow $p$-subgroup of $G$ and $Q=\langle x \rangle $ is a cyclic
 Sylow $q$-subgroup of $G$. Moreover,  $F_{q}(G)=G$ and
 $F_{p}(G)=P\langle x^{q} \rangle$.}

{\bf Lemma 3.2.} {\sl Every $ {\frak{N}}_{\sigma}$-critical group is a Schmidt group. }

{\bf Proof.} Let $G$ be an $ {\frak{N}}_{\sigma}$-critical group.
  It is clear that  for some $i$, $G$ is an    ${\frak{N}}_{\sigma _{0}}$-critical 
group, where  $\sigma _{0}=\{\sigma_{i}, \sigma_{i}'\}$. Hence  $G$ is a Schmidt group by 
\cite{bel}. The lemma is proved.

Recall that $G$ is said to be:  a \emph{$D_{\pi}$-group} if $G$ possesses a Hall 
$\pi$-subgroup $E$ and every  $\pi$-subgroup of $G$ is contained in some 
conjugate of $E$;  a \emph{$\sigma$-full group
 of Sylow type}  \cite{alg12} if every subgroup $E$ of $G$ is a $D_{\sigma _{i}}$-group for every
$\sigma _{i}\in \sigma (E)$.

{\bf Lemma 3.3 } (See Theorems A and B in \cite{2}).
   {\sl Suppose that  $G$ is $\sigma$-soluble, then $G$ is a 
$\sigma$-full group  of Sylow type and,   for every  $\Pi \subseteq \sigma$, $G$  has
  a Hall $\Pi$-subgroup $E$
 and every $\Pi$-subgroup of $G$ is contained in  conjugate of $E$. }

{\bf Lemma 3.4. } {\sl  If  $H$ is a normal 
subgroup of $G$ and  $H/H\cap \Phi (G)$ is $\Pi$-closed, then 
$E$ is $\Pi$-closed.    }

{\bf Proof.}   Let $\Phi=H\cap \Phi (G)$ and  $V/\Phi $ be a normal 
Hall $\Pi$-subgroup of $H/\Phi$.  Let $D$ be a Hall $\Pi'$-subgroup of 
$\Phi$. Then  $D$ is a normal Hall  $\Pi'$-subgroup of $V$, so 
$V$ has a Hall $\Pi$-subgroup, $E$ say,  by the Schur-Zassenhaus  
theorem.  It is clear 
that $V$ is $\pi$-soluble,  where $\pi =\bigcup _{\sigma _{i}\in \sigma (D)}\sigma _{i}$, 
 so any two Hall $\Pi$-subgroups of $V$ are conjugated in $V$. Therefore by the Frattini Argument 
we have $G=VN_{G}(E)=(E(H\cap \Phi (G)))N_{G}(E)=N_{G}(E)$.  Thus $E$   is 
normal in $G$.  The lemma is proved.

{\bf Lemma 3.5.} {\sl  Let $T/L$ be  a non-identity section of $G$. Then }

(1) {\sl $\Gamma _{H\sigma }(T/L)$ is a subgraph of $\Gamma _{H\sigma }(G)$, and}

(2) {\sl $\Gamma _{{\mathfrak{N}_{\sigma }}}(T/L)$ is a subgraph of
 $\Gamma _{{\mathfrak{N}_{\sigma }}}(G)$.}

 (3) {\sl If  $T=G$ is $\sigma$-soluble, then $\Gamma _{\sigma Hal}(G/L)$  
is a subgraph of $\Gamma _{\sigma Hal}(G)$.  }

{\bf Proof. } (1) Let $(\sigma _{i}, \sigma _{j})\in E(\Gamma _{H\sigma 
}(T))$, that is, $\sigma _{j}\in \sigma (T/F_{\{\sigma _{i}\}}(T))$. 
In view of Proposition 2.1, $F_{\{\sigma _{i}\}}(G)\cap 
T\leq F_{\{\sigma _{i}\}}(T)$, so 
from the isomorphism $T/(F_{\{\sigma _{i}\}}(G)\cap T)\simeq 
F_{\{\sigma _{i}\}}(G)T/F_{\{\sigma _{i}\}}(G)$ we get that $\sigma 
_{j}\in \sigma (TF_{\{\sigma _{i}\}}(G)/F_{\{\sigma _{i}\}}(G))$ and hence
 $(\sigma _{i}, \sigma _{j})\in E(\Gamma _{H\sigma 
}(G))$. Thus  $\Gamma _{H\sigma }(T)\subseteq \Gamma _{H\sigma }(G)$.

Now let   $(\sigma _{i}, \sigma _{j})\in E(\Gamma _{H\sigma 
}(T/L))$, that is, $\sigma _{j}\in \sigma ((T/L)/F_{\{\sigma _{i}\}}(T/L))$. 
Observe that 
 $F_{\{\sigma _{i}\}}(T)L/L\leq F_{\{\sigma _{i}\}}(T/L)$  by Proposition 2.1
 since $$F_{\{\sigma _{i}\}}(T)L/L\simeq 
F_{\{\sigma _{i}\}}(T)/(F_{\{\sigma _{i}\}}(T)\cap L)$$ and so 
$$ \sigma _{j}\in \sigma ((T/L)/(F_{\{\sigma _{i}\}}(T)L/L))=\sigma (T/F_{\{\sigma _{i}\}}(T)L)
   =$$$$=
 \sigma ((T/F_{\{\sigma_{i}\}}(T))/(F_{\{\sigma_{i}\}}(T)L/F_{\{\sigma_{i}\}}(T))) 
\subseteq  \sigma (T/F_{\{\sigma _{i}\}}(T)),$$  so 
 $(\sigma _{i}, \sigma _{j})\in E(\Gamma _{H\sigma 
}(T))$. Thus  $$\Gamma _{H\sigma }(T/L)\subseteq \Gamma _{H\sigma }(T)\subseteq 
 \Gamma _{H\sigma }(G).$$

(2)  It is clear that  $ \Gamma _{{\mathfrak{N}_{\sigma }}}(T) \subseteq 
 \Gamma _{{\mathfrak{N}_{\sigma }}}(G).$ Now let $(\sigma _{i}, \sigma 
_{j})\in E(\Gamma _{{\mathfrak{N}_{\sigma }}}(T/L))$, that is, $T/L$ has an  
${\mathfrak{N}}_{\sigma }$-critical subgroup $H/L$ such that 
$\sigma _{i} \in \sigma (H/L)$ and  
$\sigma _{j}\in \sigma ((H/L)/F_{\{\sigma _{i}\}}(H/L))$.
From Lemmas 3.1 and 3.2 it follows that $H/L$ is a $p$-closed Schmidt 
group with $\{p, q\}=\pi (H/L)$, where $p\in \sigma _{i}$ and $q
\in \sigma _{j}$.
Let $U$ be a minimal supplement of $L$ in $H$. Then $U\cap L\leq 
\Phi (U)$, so $\pi(U)=\{p, q\}$ and 
$U$ is a $p$-closed non-nilpotent  group by Lemma 3.4. Moreover, $U$ is not $\sigma$-nilpotent,
 so $U$ has an  ${\mathfrak{N}}_{\sigma }$-critical subgroup $A$.
From Lemmas 3.1 and 3.2 it follows  that $A=P\rtimes Q$  is a $p$-closed Schmidt 
group 
 with $\{p, q\}=\pi (A)$ for some $p\in \sigma _{i}$ and $q\in \sigma _{j}$, where   $P$ 
 is a Sylow $p$-subgroup of $A$ and $Q=\langle x \rangle $ is a cyclic
 Sylow $q$-subgroup of $A$ and  $F_{q}(A)=A$ and
 $F_{p}(A)=P\langle x^{q} \rangle$, which implies that $F_{\{\sigma _{i}\}}(A)=F_{p}(A) $
    and so 
 $\sigma _{j}\in \sigma (A/F_{\{\sigma _{i}\}}(A))$. Therefore $(\sigma _{i}, \sigma 
_{j})\in E(\Gamma _{{\mathfrak{N}_{\sigma }}}(T))$.  Hence 
$$\Gamma _{{\mathfrak{N}_{\sigma }}}(T/L) \subseteq \Gamma 
_{{\mathfrak{N}_{\sigma }}}(T) \subseteq \Gamma _{{\mathfrak{N}_{\sigma }}}(G).$$

(3)   Let $(\sigma _{i}, \sigma 
_{j})\in E(\Gamma _{{\mathfrak{N}_{\sigma }}}(G/L))$, that is, 
 for some Hall 
$\sigma _{i}$-subgroup $H/L$ of $G/L$ we have  $\sigma _{j}\in \sigma 
(N_{G/L}(H/L)/(H/L)C_{G/L}(H/L))$.  In view of  Lemma 3.3, for some Hall 
$\sigma _{i}$-subgroup $U$ of $G$ we have $H/L=UL/L$.  Moreover, from 
Lemma 3.3 it follows that $N_{G/L}(H/L)=N_{G}(U)L/L$. It is clear also 
that $C_{G}(U)L/L\leq C_{G/L}(H/L)$. 
Then for some $\sigma _{j}$-element $xL=gL$ of $N_{G/L}(H/L)$, where 
$g $ is a $\sigma _{j}$-element of $\in N_{G}(U)$, we have $g\not \in C_{G}(U)$.  
Hence $g\not \in C_{G}(U)U$  since $C_{G}(U)$ is normal in $C_{G}(U)U$  
and  $U$ is a $\sigma _{i}$-group, where  $i\ne j$. But then 
 $\sigma _{j}\in \sigma 
(N_{G}(U)/UC_{G}(U))$, so $(\sigma _{i}, \sigma 
_{j})\in E(\Gamma _{{\mathfrak{N}_{\sigma }}}(G))$. 
Therefore $$\Gamma _{\sigma Hal}(G/L)\subseteq \Gamma _{\sigma Hal}(G).$$

The lemma is proved.

 {\bf Proof of Proposition 1.2.}  (1)   
Assume that   $(\sigma _{i}, \sigma _{j})\in E(\Gamma _{{\sigma }Hal}(G))$ and 
let $H$ be a Hall   $\sigma _{i}$-subgroup of $G$ such that 
  $\sigma _{j}\in \sigma (N_{G}(H)/HC_{G}(H))$.  Then $i\ne j$ and for 
some non-identity 
 $\sigma _{j}$-element  $gHC_{G}(H)$  of  
$N_{G}(H)/HC_{G}(H)$, there is a $\sigma _{j}$-element  $a\ne 1$ of  
$N_{G}(H)$ such that $gHC_{G}(H) =aHC_{G}(H)$.  Hence  $a\not \in C_{G}(H)$, 
so  $E=H\rtimes  \langle a\rangle$ is non-$\sigma$-nilpotent 
$\sigma _{i}$-closed group.   
 Let $A$ be an ${\mathfrak{N}}_{\sigma 
}$-critical subgroup of $E$. Then, in view of Lemmas 3.1 and 3.2, 
$A=P\rtimes Q$, where $P$ is a Sylow $p$-subgroup and $Q$ is a Sylow $q$-subgroup of $A$ 
 for some  $p\in \sigma _{i}$ and   $q\in \sigma _{j}$ and $Q\nleq F_{\{\sigma _{i}\}}(A)$.
 Then    $\sigma _{j}\in \sigma (A/F_{\{\sigma _{i}\}}(A))$.
Hence  $(\sigma _{i}, \sigma _{j})\in E(\Gamma _{{\mathfrak{N}_{\sigma }}}(G))$. 
Thus we have   $\Gamma _{{\sigma }Hal}(G)  \subseteq
 \Gamma _{{\mathfrak{N}}_{\sigma }}(G)$.

Assume that $(\sigma _{i}, \sigma _{j})\in E(\Gamma 
_{{\mathfrak{N}}_{\sigma }}(G))$. Then, in view of Lemmas 3.1 and 3.2,
 $G$ 
has a Schmidt subgroup $A=A_{p}\rtimes A_{q}$ such that $p\in \sigma _{i}$, 
$q\in \sigma _{j}$ ($i\ne j$) and  $A_{q}\nleq  
F_{\{\sigma _{i}\}}(A)$. Therefore  $A_{q}\nleq  
F_{\{\sigma _{i}\}}(G)\cap A\leq  F_{\{\sigma _{i}\}}(A)$. Hence $(\sigma _{i}, 
\sigma _{j})\in E(\Gamma _{H\sigma }(G)$, so $\Gamma _{{\mathfrak{N}}_{\sigma }}(G)\subseteq 
  \Gamma _{H\sigma }(G)$.

Now we show that  both inclusions   in (*) may be strict. Let $ p < q$ be primes such that
 $p$ divides $q-1$ and  let $C_{q}\rtimes C_{p}$ be   a non-abelian group of 
order $pq$. Let $G=C_{p}\wr (C_{q}\rtimes C_{p})=K\rtimes (C_{q}\rtimes 
C_{p})$, where $K$ is the base group of the regular wreath product $G$. 
Let $\sigma =\{\sigma _{1}, \sigma _{2}\}$, where $\sigma _{1}=\{p\}$ and $\sigma _{2}=
 \{p\}'$. Then $K=F_{\{\sigma _{1}\}}(G)$, so $(\sigma _{1}, \sigma _{1})\in 
E(\Gamma _{H\sigma 
}(G))$. Hence   $\Gamma _{{\mathfrak{N}}_{\sigma }}(G) \subset \Gamma _{H\sigma }(G)$.
 We have also that $(\sigma _{1}, \sigma _{2})\in 
E(\Gamma _{H\sigma }(G))$ but $(\sigma _{1}, \sigma _{2})\not \in 
E(\Gamma _{{\sigma 
}Hal}(G))$ since $C_{q}\nleq N_{G}(KC_{p})$. Hence $\Gamma _{{\sigma 
}Hal}(G) \subset  \Gamma _{{\mathfrak{N}}_{\sigma }}(G).$

 Assume that $\Gamma _{{\sigma 
}Hal}(G) = \Gamma _{H\sigma 
}(G)$. Then  $\Gamma _{H\sigma 
}(G)$ has no loops in every vertex $\sigma _{i}\in \sigma (G)$. It follows that 
$G/F_{\{\sigma _{i}\}}(G)$ is a $\sigma _{i}'$-group. On the other hand, by Proposition 2.3,   
$F_{\{\sigma _{i}\}}(G)$ is $\sigma _{i}'$-closed and so $G$ is $\sigma$-soluble
 with  $l_{\sigma _{i}}(G)\leq 1$ for all $\sigma _{i}\in \sigma (G)$.

 Finally, we show that if   $G$ is $\sigma$-soluble  and $l_{\sigma 
_{i}}(G)\leq 1$ for all $\sigma _{i}\in \sigma (G)$, then  $\Gamma _{{\sigma 
}Hal}(G) = \Gamma _{H\sigma }(G)$.  The inclusion  $\Gamma _{{\sigma 
}Hal}(G) \subseteq  \Gamma _{H\sigma }(G)$ follows from Part (1). 
Next we show that $\Gamma _{H\sigma }(G)\subseteq \Gamma _{{\sigma 
}Hal}(G)$. 
If $|\sigma (G)| =1$, it is evident.
 Now assume that  
 $\sigma (G)=\{\sigma _{1}, \ldots ,\sigma _{t}\}$, where $t > 1$, and let 
$(\sigma _{i}, \sigma _{j})\in E(\Gamma _{H\sigma 
}(G))$.   
 Since $l_{\sigma 
_{i}}(G)\leq 1$, a Hall $\sigma _{i}$-subgroup $H$ of   $F_{\{\sigma 
_{i}\}}(G)$ is  also a Hall $\sigma _{i}$-subgroup of $G$. 
 Moreover, 
every two  Hall $\sigma _{i}$-subgroups  of $F_{\{\sigma 
_{i}\}}(G)$ are conjugate in $F_{\{\sigma
_{i}\}}(G)$ by Lemma 3.3 and  so $G= F_{\{\sigma _{i}\}}(G) N_{G}(H)$ by the Frattini 
Argument, which implies that $\sigma _{j}\in \sigma (N_{G}(H)/(N_{G}(H)\cap F_{\{\sigma _{i}\}}(G)))$.    Hence  for some non-identity  $\sigma _{j}$-subgroup $U$ of $G$ 
we have $U\nleq O_{\sigma _{i}'}(G)\leq F_{\{\sigma _{i}\}}(G) $ and $U\leq N_{G}(H)$. Then 
$UO_{\sigma _{i}'}(G)/O_{\sigma _{i}'}(G)$ is a non-identity group, which 
implies that
 $$UO_{\sigma _{i}'}(G)/O_{\sigma _{i}'}(G)\nleq
 C_{G/O_{\sigma _{i}'}(G)}(F_{\{\sigma _{i}\}}(G)/O_{\sigma _{i}'}(G))=C_{G/O_{\sigma _{i}'}(G)}(HO_{\sigma _{i}'}(G)
/O_{\sigma _{i}'}(G)).$$ 
by Theorem 3.2 in \cite[Chapter 6]{Gor}. But then $U\nleq C_{G}(H)$, so  $U\nleq HC_{G}(H)$ 
since $i\ne j$. 
 Hence  $\sigma _{j}\in \sigma (N_{G}(H)/HC_{G}(H))$, which implies that 
$(\sigma _{i}, \sigma _{j})\in E(\Gamma _{{\sigma 
}Hal}(G))$. Therefore  $\Gamma _{H\sigma 
}(G)\subseteq  \Gamma _{{\sigma 
}Hal}(G)$, so   $\Gamma _{{\sigma 
}Hal}(G) = \Gamma _{H\sigma 
}(G)$.  Therefore the Statement (2) holds.

The proposition is proved.

{\bf Lemma 3.6.} {\sl If $G$ is a $\sigma$-soluble group, then for every
 $\sigma _{i}\in \sigma (G)$ the group $G$ possesses a maximal subgroup $M$ such that 
$|G:M|$ is a $\sigma _{i}$-number.}

{\bf Proof. }  We prove this lemma by induction on $|G|$. Let $R$ be a minimal normal 
subgroup of $G$. Then  $R$ is a  $\sigma _{j}$-group for some $j$ since 
$G$ is  $\sigma$-soluble.

For every $\sigma _{i}\in 
\sigma (G/R) $ the group 
  $G/R$ possesses a maximal subgroup $M/R$ such that $|G/R:M/R|=|G:M|$ is a $\sigma 
_{i}$-number. Therefore  the assertion holds for every $\sigma _{i}\in 
\sigma (G)\setminus \{\sigma _{j}\}$. If  $R\leq \Phi (G)$, then 
 $R$ is not a Hall subgroup of $G$, so 
 $\sigma (G) = \sigma (G/R)$ and hence the assertion holds for $G$. 
Finally, if $R\nleq \Phi(G)$, then for some maximal subgroup $M$ of $G$ we 
have $G=RM$ and so  $|G:M|$ is a $\sigma _{j}$-number. 
The lemma is proved.

A set  $ {\cal H}$ of subgroups of $G$ is said to be a  \emph{ complete
Hall $\sigma $-set} of $G$ \cite{alg12} if   every  member $\ne 1$  of  ${\cal H}$ is a
 Hall $\sigma _{i}$-subgroup
 of $G$ for some $i$ and ${\cal H}$ contains exactly one Hall  $\sigma
_{i}$-subgroup of $G$ for every $\sigma _{i} \in \sigma (G)$.

 {\bf Proof of Proposition 1.11.}  Suppose that this proposition is false and let $G$ be a 
counterexample of minimal order. Then $\Pi \cap \sigma(G) \ne \emptyset \ne \Pi'
 \cap \sigma (G)$.

By Lemma 3.3, $G$ possesses 
  a  complete Hall  $\sigma$-set $\{H_{1}, \ldots ,
H_{t} \}$.
Without loss of generality we can assume that $H_{i}$ is a  $\sigma_{i}$-group for
 all $i=1, \ldots , t$ and that $\Pi \cap \sigma (G) =\{\sigma_{1}, \ldots , \sigma_{n}\}$ and 
$\Pi'  \cap \sigma (G)=\{\sigma_{n+1}, \ldots , \sigma_{t}\}$.   Let  $R$ be 
 a minimal normal subgroup of   $G$. Then $R$ is a  $\sigma_{j}$-group for 
some $j$ since $G$ is $\sigma$-soluble by hypothesis.

(1) {\sl  If either $R\leq \Phi (G)$ or $\sigma_{j}\in \Pi$, then  $G/R$ is a  $\Pi'$-closed  
Schmidt group. }

The choice of $G$ and Lemma 3.4 shows that  $G/R$ is not  $\Pi$-closed.
  On the other hand, every 
maximal subgroup $M/R$ of $G/R$ is   $\Pi$-closed since $M$ is $\Pi$-closed by hypothesis.
 Hence the hypothesis holds for $G/R$. Therefore the 
choice of $G$   implies that   $G/R$ is a  $\Pi'$-closed  
Schmidt group.

(2)  {\sl  $\Phi (G)=1 $,  $R$ is  the unique 
minimal normal subgroup of $G$ and $R$ is a $\Pi'$-group.  Hence $C_{G}(R)\leq R$. }

Suppose that    $R\leq \Phi (G)$.  Then  $G/R$ is a  $\Pi'$-closed  
Schmidt group by Claim (1), so in view  of  Lemmas 3.1, 3.2 and 3.4,  $G= 
H_{1}\rtimes H_{2}=P\rtimes Q$, where  $H_{1}=P$ is a $p$-group and $H_{2}=Q$
 is a $q$-group for some primes $p\in \Pi'$ and $q\in \Pi$. Therefore, in fact, $G$  is 
not $p$-nilpotent but 
 every maximal subgroup of $G$ is 
$p$-nilpotent.  Hence  $G$ is a $p$-closed  
Schmidt group by \cite[IV, Satz 5.4]{hupp}, a contradiction.  Therefore $R\nleq \Phi (G)$.

Now  assume that $G$ has a minimal normal subgroup $L\ne R$.
Then there are maximal subgroups $M$ and $T$   of $G$ such that $LM=G$ and 
$RT=G$. By hypothesis, $M$ and $T$ are  $\Pi$-closed. Hence 
$G/L\simeq LM/L\simeq M/(M\cap L)$ is $\Pi$-closed. Similarly, 
$G/R$ is $\Pi$-closed and so $G\simeq G/(L\cap R)$ is $\Pi$-closed, a contradiction.
  Hence $R$ is  the unique 
minimal normal subgroup of $G$. Hence  $C_{G}(R)\leq R$ by \cite[Chapter A, Theorem 15.6]{DH}.  It is also clear  that  $R$ is a 
$\Pi'$-group.

(3)  {\sl   $|\sigma (G)|=2$. Hence $G=H_{1}H_{2}$ and $n=1$.}

It is clear that $|\sigma (G)| > 1$. Suppose that $|\sigma (G)| > 2 $.
Then, since $G$ is $\sigma$-soluble, there are maximal subgroups $M_{1}$, 
$M_{2}$ and $M_{3}$ whose indices $|G:M_{1}|$, $|G:M_{2}|$  and 
$|G:M_{3}|$   are pairwise  $\sigma$-coprime by Lemma 3.6.

Hence $G=M_{1}M_{2}=M_{2}M_{3}=M_{1}M_{3}$ and for some $i$ and $j$, say $i=1$ and $j=2$, we have 
$R\leq M_{1}\cap M_{2}$. Then   $O_{\Pi}(M_{1})=1=O_{\Pi}(M_{2})$ by Claim (2).
  But by hypothesis, $M_{1}$  and $ M_{2}$  are 
 $\Pi$-closed and hence  $M_{1}$  and $ M_{2}$ are  $\Pi'$-groups. 
$G=M_{1}M_{2}$ is  a  $\Pi'$-group.  This 
contradiction completes the proof of Claim (3).

  {\sl  Final contradiction.}  In view of Claims (2) and (3), 
$C_{G}(R) \leq R \leq H_{2}$. The subgroup $RL$ is $\Pi$-closed 
 for  every proper subgroup $L$ of $H_{1}$. Hence $L=1$, so 
 $|H_1{}|$ is a prime and 
   $RH_{1}=G$ since $R\leq H_{2}$ and every proper subgroup of $G$
 is $\Pi$-closed. Therefore $R=H_{2}$, so $R$ is not abelian since  $G$ is  not  
a $\Pi'$-closed  Schmidt group.    From the Frattini Argument it follows  that  for any prime 
$p$ dividing $|R|$ there is a Sylow $p$-subgroup $P$ of $G$ such that 
$H_{1}\leq N_{G}(P)$, so  $PH_{1}=H_{1}P$.
But $H_{1}P < G$, so   $H_{1}P=H_{2}\times  P$. Therefore $R\leq 
N_{G}(H_{1})$ and hence $G=R\times H_{1} =H_{1}\times H_{2}
$ is $\sigma$-nilpotent. This final contradiction completes the proof of the result.

{\bf Proof of Proposition 1.9. }  Assume that this theorem is false and let $G$ be a 
counterexample of minimal order.  Then $\Pi _{1}\ne \emptyset\ne \Pi_{2}$.

First suppose that  $\Gamma=\Gamma _{H\sigma }(G)$ and let  
$\sigma _{j}\in \Pi_{2}$. The   $G/F_{\{\sigma _{j}\}}(G)$ is a $\Pi_{2}$-group  by hypothesis. 
 On the other hand,  $\Pi _{1}\ne \emptyset$. 
 Hence  $F_{\{\sigma _{j}\}}(G) \ne 1$. 

 Next show that  $F_{\{\sigma _{j}\}}(G) =G$. Indeed, assume that $F_{\{\sigma _{j}\}}(G)  < G$.
In view of Lemma 3.4(1),  
$\Gamma _{H\sigma }(F_{\{\sigma _{j}\}}(G))\subseteq \Gamma _{H\sigma }(G)$
 and so the hypothesis holds for
 $F_{\{\sigma _{j}\}}(G)\ne 1$.
Therefore $F_{\{\sigma _{j}\}}(G)$ possesses a normal Hall $\Pi_{1}$-subgroup $H$ by the 
choice of $G$. Then $H$ is characteristic in $F_{\{\sigma _{j}\}}(G)$, so it is normal in 
$G$. 
Finally, $H$ is a Hall $\Pi_{1}$-subgroup of $G$ since 
 $G/F_{\{\sigma _{j}\}}(G)$ is a $\Pi_{2}$-group  and so   $G$ is $\Pi_{1}$-closed, a contradiction.
 Therefore  $F_{\{\sigma _{j}\}}(G)=G$ for all 
$\sigma _{j}\in \Pi_{2}$. 
 Then, by Proposition 2.3,   $G$ has a 
normal $\sigma _{j}$-complement $V_{j}$ for all $\sigma _{j}\in \Pi_{2}$.
  But then $\bigcap _{\sigma _{j}\in \Pi_{2}}V_{j}$ 
 is a normal Hall $\Pi _{1}$-subgroup of $G$, so $G$  is $\Pi _{1}$-closed, a   contradiction. 
Thus this theorem in the case when  $\Gamma=\Gamma _{H\sigma }(G)$ is 
true.

Now suppose that $G$ is $\sigma$-soluble and   $\Gamma=\Gamma _{{\mathfrak{N}}_{\sigma 
}}(G)$. The hypothesis holds for every subgroup of $G$, so the choice of $G$ implies that 
$G$ is a minimal non-$\Pi_{1}$-closed group. Then $G$ is a $\Pi_{2}$-closed 
Schmidt subgroup  of $G$ by Proposition 1.11. But then, in view of Lemmas 
3.1 and 3.2,  for some $\sigma _{j}\in \Pi_{2}$ and $\sigma _{i}\in \Pi _{1}$ we 
have  $(\sigma _{i}, \sigma _{j})\in E(\Gamma)$, a contradiction.
Thus this theorem in the case when  $\Gamma=\Gamma _{{\mathfrak{N}}_{\sigma 
}}(G)$ is also true.

The proposition is proved.

\section{Proofs of Theorem 1.4, 1.7, and 1.13}

We say, following \cite{1}, that  $G$  is \emph{$\sigma$-dispersive} or \emph{ $G$ possesses 
$\sigma$-Hall tower} if   $G$ 
has a  normal  series   $1= G_{1} <  G_{2}  < \cdots   <  G_{t-1} < 
G_{t}=G$   and  a
 complete Hall $\sigma$-set ${\cal H}=\{H_{1}, \ldots ,
H_{t} \}$   such that    $G_{i}H_{i}=G_{i+1}$ for all $i=1,  \ldots , t-1$.

The following lemma is a corollary of Lemma 3.3.  

{\bf Lemma 4.1.} {\sl $G$ is $\sigma$-dispersive if and only if $G$ satisfies the
 Condition (1) in Theorem 1.7.}

{\bf Proof of Theorem 1.4. }    (1) $\Rightarrow$ (2) 
  By hypothesis and Lemma 4.1, $G$ 
has a  normal  series   $1= G_{1} <  G_{2}  < \cdots   <  G_{t-1} < 
G_{t}=G$   and  a
 complete Hall $\sigma$-set ${\cal H}=\{H_{1}, \ldots ,
H_{t} \}$   such that    $G_{i}H_{i}=G_{i+1}$ for all $i=1,  \ldots , t-1$. 
 We can assume without loss of generality that $H_{i}$ is $\sigma _{i}$-group 
for all $i$. From Propositions 2.1 an 2.3 it follows that  
$G_{i+1}=F_{\{\sigma _{i}\}}(G)$, so $\sigma _{i}\not \in
 (G/F_{\{\sigma_{i}\}}(G))$ since $F_{\{\sigma _{i}\}}(G))$
 contains a Hall $\sigma _{i}$-subgroup $H_{i}$ of $G$.
Thus  $\Gamma _{H\sigma }(G)$ has no circuits.

 (2) $\Rightarrow$ (3) By hypothesis,  $\Gamma _{H\sigma }(G)$ has no 
loops, so $G$ is $\sigma$-soluble and
 $\Gamma _{{\mathfrak{N}_{\sigma }}}(G)=\Gamma _{H\sigma }(G)$ by 
Proposition 1.2(2) 
 and hence $\Gamma _{{\mathfrak{N}_{\sigma }}}(G)$ has no circuits.

 (3) $\Rightarrow$ (1)  We prove this implication by induction on $|G|$.  
 First note that since $\Gamma _{{\mathfrak{N}}_{\sigma }}(G)$  has no 
circuits, there exists $\sigma _{i}\in \sigma (G)$  such that for every 
$\sigma _{j}\in \sigma (G)\setminus \{\sigma _{i}\}$ the group $G$ has no 
an   ${\mathfrak{N}_{\sigma }}$-critical subgroup $A$ such that $\sigma _{j}\in \sigma
 (A/F_{\{\sigma _{i}\}}(A))$. 
 Taking now $\Pi =\{\sigma _{i}\}'$ we have 
$\Pi'=\{\sigma _{i}\}$ and so $G$ has a normal Hall  $\Pi$-subgroup $V$ by 
Proposition 1.9.  The hypothesis holds for $V$ by Proposition 1.2, so $V$ is $\sigma$-dispersive
 by 
induction and Lemma 4.1. Moreover, every normal Hall subgroup of $V$ is characteristic 
in $V$ and so normal in $G$. It follows that $G$ is $\sigma$-dispersive, so (3) $\Rightarrow$ (1).

The theorem is proved.

{\bf Lemma 4.2} (See   Proposition 1.1 in \cite{2}). Suppose that $G = AB = 
Bc = AC$, where $A$, $B$ and $C$ are  $\sigma$-soluble subgroups of G. If the
 three indices $|G:A|$, $|G:B|$ and $|G:C|$ are pairwise $\sigma$-coprime, then G is $\sigma$--soluble.

{\bf Lemma 4.3.} {\sl Suppose that $G=AB=AC=BC$ is
 $\sigma$-soluble, where $A$, $B$ and $C$ are subgroups of $G$. Then for every
 $\Pi \subseteq \sigma$, there exists an element $x\in G$ and Hall   $\Pi$-subgroups 
$A_{\Pi}$, $B_{\Pi}$,  $C_{\Pi}$ and $G_{\Pi}$ of $A$, $B$, $C^{x}$ and 
$G$, respectively, such that $G_{\Pi}=A_{\Pi}B_{\Pi}=B_{\Pi}C_{\Pi}=A_{\Pi}C_{\Pi}$.  }

{\bf Proof.}  See the proof of Lemma 1(3) in \cite{pen} or Lemma 5 in  \cite{555}.

  We write $O_{\Pi}(G)$ to denote the largest normal $\Pi$-subgroup of $G$.

{\bf Lemma 4.4.} {\sl  If $H$ is a subnormal $\Pi$-subgroup of $G$, then $H\leq  O_{\Pi}(G)$.}

{\bf Proof.} Let $H=H_{0} \trianglelefteq H_{1} \trianglelefteq \cdots 
\trianglelefteq 
H_{t-1} \trianglelefteq H_{t}=G$. We prove this lemma by induction on $t$. 
First observe that $O_{\Pi}(H_{t-1})$ is characteristic in $H_{t-1}$, so 
it is normal in $G$, so  $O_{\Pi}(H_{t-1})\leq O_{\Pi}(G)$. On the other 
hand, 
 we have 
 $H\leq O_{\Pi}(H_{t-1})$ by induction. The lemma is proved.

{\bf Lemma 4.5.} {\sl  If  $G=AB=BC=AC$ is $\sigma$-soluble, where the subgroups $A$, $B$ and 
$C$ are $\Pi$-closed, then $G$ is also  $\Pi$-closed.}

{\bf Proof.} Assume that this lemma is false and let $G$ be a 
counterexample of minimal order. Then $\Pi\cap \sigma (G)\ne \emptyset \ne 
\Pi\cap \sigma (G)$. Then for at least two of the subgroups $A$, $B$, $C$, for
 $A$ and $B$ say, we have  $\Pi\cap \sigma (A)\ne \emptyset \ne 
\Pi\cap \sigma (B)$. Let $A_{\Pi}$ and $B_{\Pi}$ be the Hall $\Pi$-subgroups of $A$ and $B$, 
respectively.

(1) {\sl If $R$ is a minimal  normal subgroup of $G$, then $R$  is 
 the unique minimal normal subgroup of $G$, $R\nleq \Phi 
(G)$ and $\sigma (R)\not \subseteq  \Pi$}
 (Since the hypothesis holds for $G/R$, this follows from the choice of $G$).

(2) {\sl   $A_{\Pi}B_{\Pi}^{x}=B_{\Pi}^{x}A_{\Pi}$ is 
a Hall $\Pi$-subgroup  of $G$  for all $x\in G$. }

From $G=AB$ it follows that for some Hall $\Pi$-subgroup $G_{\Pi}$ of $G$ we
 have   $G_{\Pi}=A_{\Pi}B_{\Pi}$ by \cite[Theorem 1.1.19]{prod}  and Lemma 3.3.
Moreover, from  $G=AB$ we get also that $G=AB^{x}$ for all $x\in G$ by \cite[p. 675]{hupp}. 
Then  $A_{\Pi}B_{\Pi}^{x}=B_{\Pi}^{x}A_{\Pi}$ is 
a Hall $\Pi$-subgroup  of $G$ since $B^{x}$ is also $\Pi$-closed.
 Therefore $A_{\Pi}$ permutes with  
$B_{\Pi}^{x}$ for all $x\in G$.

(3) {\sl $[A_{\Pi}, B_{\Pi}^{x}]=1$ for all $x\in G$.}

Let $L=B_{\Pi}^{x}$. Then  $A_{\Pi}$ 
permutes with $L^{z}=(B_{\Pi}^{x})^{z}=B_{\Pi}^{xz}$ for all $z\in G$ by Claim (2). 
Let    $H=(A_{\Pi})^{L}\cap L^{A_{\Pi}}$.  Then $H$ is  subnormal subgroup 
of $G$ by  \cite[Lemma 1.1.9(2)]{prod} and Claim (2).  But  $H\leq A_{\Pi}L=LA_{\Pi}$,
 so $H$ is a $\Pi$-group and hence in the case when $H\ne 1$ we have 
 $O_{\Pi}(G)\ne 1$ by Lemma 4.4. But then  $R\cap  O_{\Pi}(G)\ne 1$ by Claim (1) and so 
$R$ is a $\Pi$-group,  contrary to Claim (1).  Therefore $H=1$, so we have 
  $[A_{\Pi}, L]\leq [(A_{\Pi})^{L}, L^{A_{\Pi}}]\leq H=1$. Hence (3) holds.

 {\sl Final contradiction.} Fro
m Claim  (3) it follows that $[A_{\Pi}, 
(B_{\Pi})^{G}]=1$, so $(B_{\Pi})^{G}\leq C_{G}(A_{\Pi})$. But $B_{\Pi}\ne 1$  
since $ \Pi\cap \sigma (B)\ne \emptyset$. Hence from Claim (1) it follows 
that $R\leq (B_{\Pi})^{G}\leq  C_{G}(A_{\Pi})$. Hence  $1  < A_{\Pi}\leq  C_{G}(R)\leq R$
by Claim (1) and \cite[Chapter A, Theorem 15.6]{DH}. But this is impossible since 
$\sigma (R)\not \subseteq  \Pi$ by Claim (1). This contradiction completes 
the proof of the result.

{\bf Proof of Theorem 1.7.}   Let $\Gamma=\Gamma _{{\mathfrak{N}_{\sigma 
}}}(G)$ and let $\Gamma _{1}=\Gamma _{{\mathfrak{N}_{\sigma }}}(A)
\bigcup\Gamma _{{\mathfrak{N}_{\sigma }}}(B) 
 \bigcup 
\Gamma _{{\mathfrak{N}_{\sigma }}}(C)$ and $\Gamma _{2}=\Gamma _{H\sigma }(A) \bigcup\Gamma _{H\sigma }(B)
  \bigcup \Gamma _{H\sigma }(C)$.

(1) In view of Lemma 3.4(2), it is enough to show that  
$\Gamma  \subseteq  \Gamma _{1}.$
Assume that this is false and let $G$ be a counterexample of minimal 
order.

First we show that $|\sigma (G)|=2$. Assume that $|\sigma (G)|  > 2$ and let
$(\sigma _{i}, \sigma _{j}) \in E(\Gamma)$. Let $\Pi =\{\sigma _{i}, \sigma 
_{j}\}$. Then, in view of Lemmas  3.1 and 3.2, for some Schmidt subgroup $A$ of $G$ we have 
and $\sigma (A)=\Pi$ and 
$\sigma _{j}\in \sigma (A/F_{\{\sigma _{i}\}})$.     
Since $G$ is $\sigma$-soluble 
by hypothesis, from Lemmas 3.3 it follows that $G$ has a Hall 
$\Pi$-subgroup  $P$ such that $A\leq P$. Moreover, in view of 
Lemma 4.3, it follows that there exists an element $x\in G$ and Hall   $\Pi$-subgroups 
$A_{\Pi}$, $B_{\Pi}$,  $C_{\Pi}$ and $G_{\Pi}$ of $A$, $B$, $C^{x}$ and 
$G$, respectively, such that $G_{\Pi}=A_{\Pi}B_{\Pi}=B_{\Pi}C_{\Pi}=A_{\Pi}C_{\Pi}$. 
Then the hypothesis holds for $G_{\Pi}$ by Lemma 3.4(2) and $|G_{\Pi}| < |G|$ since $|\sigma
 (G)|  > 2$. 
Hence  the choice of $G$ implies that $$E(\Gamma _{{\mathfrak{N}_{\sigma }}}(G_{\Pi}))\subseteq 
     E(\Gamma _{{\mathfrak{N}_{\sigma }}}(A_{\Pi}))\bigcup E(\Gamma _{{\mathfrak{N}_{\sigma }}}(B _{\Pi})) 
 \bigcup E(\Gamma _{{\mathfrak{N}_{\sigma }}}(C_{\Pi})).$$ But every two 
Hall $\Pi$-subgroups of $G$ are conjugate by Lemma 3.3, so some conjugate of $A$ is a 
subgroup of $G_{\Pi}$. Hence  $(\sigma _{i}, \sigma _{j}) \in E(\Gamma 
_{1})$.  Therefore $\Gamma  \subseteq  \Gamma _{1},$  a contradiction. 
Hence $|\sigma (G)|=2$.

 IOt follows that  for some $(\sigma _{i}, \sigma _{j}) \in E(\Gamma)$ we have 
$$(\sigma _{i}, \sigma _{j}) \not \in E(\Gamma _{{\mathfrak{N}_{\sigma }}}(A))\bigcup E(\Gamma
 _{{\mathfrak{N}_{\sigma }}}(B)) 
 \bigcup E(\Gamma _{{\mathfrak{N}_{\sigma }}}(C)).$$ Therefore the 
subgroups $A$, $B$ and $C$ are  $\sigma _{j}$-closed by Proposition 1.9 and so $G$ is    
$\sigma _{j}$-closed by Lemma 3.9. But then $G$ has no a $\sigma 
_{i}$-closed  ${\mathfrak{N}}_{\sigma}$-critical subgroup and so $(\sigma 
_{i}, \sigma _{j}) \not  \in E(\Gamma)$. This contradiction completes the 
proof of the Statement (1).

(2)  In view of Lemma 3.4(1), it is enough to show that  
$\Gamma  \subseteq  \Gamma _{2}.$
Assume that this is false and let $G$ be a counterexample of minimal 
order.  Then

Observe that $G$ is $\sigma$-soluble by Lemma 4.2.   Let
$(\sigma _{i}, \sigma _{j}) \in E(\Gamma)$, that is, $\sigma _{j}\in  \sigma
 (G/F_{\{\sigma \}}(G)$. `	

First assume that $D=O_{\sigma _{i}'}(G)\ne 1$. From Proposition 2.3 it 
follows that
 $F_{\{\sigma _{i}\}}(G/D)=F_{\{\sigma _{i}\}}(G)/D =O_{\sigma 
_{i}}(G/D). $  Therefore $\sigma _{j}\in  \sigma
 ((G/D)/F_{\{\sigma \}}(G/D))$. On the other hand, $AD/D$, $BD/D$ and $CD/D$ are $\sigma$-soluble 
subgroups of $G/D$ and the  indices $|G/D:AD/D|, |G/D:BD/D|,
 |G/D:CD/D|$ are pair   $\sigma$-coprime, so the hypothesis holds for 
$G/D$. Therefore, in view of Lemma 3.4(1), 
 $$\Gamma _{H\sigma }(G/D)  \subseteq \Gamma _{H\sigma }(AD/D) 
\bigcup\Gamma _{H\sigma }(BD/D))
  \bigcup \Gamma _{H\sigma }(CD/D)$$$$=\Gamma
 _{H\sigma }(A/(D\cap A)) \bigcup\Gamma _{H\sigma }(B/(D\cap B))
  \bigcup \Gamma _{H\sigma }(C/(D\cap C))$$$$\subseteq \Gamma _{H\sigma }(A) \bigcup\Gamma _{H\sigma }(B)
  \bigcup \Gamma _{H\sigma }(C)$$  by the choice of $G$. Hence	  $\Gamma  \subseteq
  \Gamma _{2}$, a 
contradiction.  Therefore $D=1$, so $F_{\{\sigma _{i}\}}(G)=O_{\sigma 
_{i}}(G)$. 

From Lemma 3.3 and the hypothesis it follows that $G$ has a Hall $\sigma 
_{j}$-subgroup $H$ such that  for at least one of the subgroups $A$, $B$
 or $C$, for $A$ say,  we have  $H\leq A$ and $O_{\sigma 
_{i}}(G)\leq A$.  On the other hand, $C_{G}(O_{\sigma_{i}}(G))\leq  
O_{\sigma_{i}}(G)$ by Theorem 3.2 in \cite[Chapter 6]{Gor}.  Hence 
$F_{\{\sigma _{i}\}}(A)=O_{\sigma _{i}}(A)$. 

Assume that $j\ne i$. Then  $\sigma _{j}\in  \sigma
 (A/F_{\{\sigma _{i}\}}(A)=\sigma
 (A/O_{\sigma _{i}}(A)$ since $H\ne 1$, so $(\sigma _{i}, \sigma _{j}) \in E(\Gamma _{2})$. 	

Finally, assume that $j=i$. In this case we may assume without loss of 
generality that $H, O_{\sigma_{i}}(G)\leq A\cap B$. Suppose that 
$(\sigma _{i}, \sigma _{j})=(\sigma _{i}, \sigma _{i}) \not \in E(\Gamma _{2})$. Then
 $H= O_{\sigma _{i}}(A)=O_{\sigma _{i}}(B)$, so $H$ is normal in $G=AB$. But then 
$H=F_{\{\sigma _{i}\}}(G)$ and so $\sigma _{i}=\sigma _{i}\not \in \sigma (G/F_{\{\sigma _{i}\}}(G))$.
This contradiction shows that $(\sigma _{i}, \sigma _{j}) \in E(\Gamma _{2})$ in any case.  
Therefore  $\Gamma  \subseteq  \Gamma _{2}.$ Hence the Statement (2) 
holds.

The theorem is proved.

{\bf Proof of Theorem 1.12.}   The implication  (1) $\Rightarrow$  (2) is evident. The 
implication (2) $\Rightarrow$  (3)  is a corollary of Proposition 1.2.

(3) $\Rightarrow$  (4)   It is enough to show that the Condition (3) 
implies that $G$ is $\sigma$-soluble. Assume that this is false. Then $G$ 
is not $\sigma$-nilpotent, so it has a   Schmidt subgroup  $A=P\rtimes Q$, where 
   $P=A^{\frak{N}}$ 
 is a Sylow $p$-subgroup of $A$ and $Q=\langle x \rangle $ is a cyclic
 Sylow $q$-subgroup of $A$ for some $p\in \sigma _{i}$ and $q\in \sigma _{j}$ 
by Lemmas 3.1 and 3.2.
 Moreover,  $F_{q}(G)=G$ and
 $F_{p}(G)=F_{\{\sigma_{i}\}}(G)=P\langle x^{q} \rangle$. Hence $\sigma _{j}\in \sigma (A/F_{\{\sigma _{i}\}}(A))$,
 so $(\sigma _{i}, \sigma _{j})\in 
E(\Gamma _{{\mathfrak{N}}_{\sigma }}(G))=\emptyset, $  a contradiction. 
Hence the implication holds.

(4) $\Rightarrow$  (1)  Assume that this implication is false and let $G$ 
be a counterexample of minimal order.  In view of Lemma 3.5(3), 
the hypothesis holds on $G/R$ for every minimal normal subgroup $R$ of 
$G$. Therefore the choice of $G$ implies that $G/R$ is $\sigma$-nilpotent.
If $G$ has a minimal normal subgroup $N\ne R$, then $G\simeq G/(R\cap N)$ 
is $\sigma$-nilpotent by Proposition 2.1, contrary to our assumption on $G$.
Hence $R$ is the unique minimal normal subgroup of $G$ and, by Lemma 3.4, 
$\Phi(G)=1$. Therefore $C_{G}(R)\leq R$ by \cite[Chapter A, Theorem 15.2]{DH}.
Since $G$ is $\sigma$-soluble by hypothesis, $R$ is a $\sigma _{i}$-group for some 
$i$ and $G$ has a Hall $\sigma _{i}$-subgroup $H$ by Lemma 3.3. Then $R\leq H$ and $H/R$ 
is a normal Hall $\sigma _{i}$-subgroup of $G/R$ by Proposition 2.3. Therefore $N_{G}(H)=G$
and $ C_{G}(H)\leq C_{G}(R)\leq R$. But then for some $j\ne i$ we have 
$\sigma _{j}\in \sigma (N_{G}(H)/HC_{G}(H))$, that is,  $\sigma _{j}\in
 E(\Gamma _{{\sigma }Hal}(G))=\emptyset$.
 This contradiction completes the proof of the implication.

The theorem is proved.

\end{document}